\documentclass[11pt]{amsart} 

\usepackage{verbatim, latexsym, amssymb, amsmath,color}
\usepackage{epsfig}
\def\R{\mathbb R}
\def\RP{\mathbb {RP}}
\def\N{\mathbb N}
\def\Z{\mathbb Z}

\def\vol{\mathrm{vol}}

\newtheorem{thm}{Theorem}[section]

\newtheorem*{coro}{Corollary} 
\newtheorem*{Wconj}{Willmore Conjecture (1965)}
\newtheorem*{lconj}{Freedman--He--Wang Conjecture (1994)}
\newtheorem*{yconj}{Yau's Conjecture (1982)}
\theoremstyle{definition}

\newtheorem{remark}[thm]{Remark}

\title{New applications of  Min--max Theory}
\author{Andr\'e Neves}
\address{Imperial College \\ Huxley Building \\ 180 Queen's Gate \\ London SW7 2RH \\ United Kingdom}
\email{a.neves@imperial.ac.uk}
\thanks{The author author was partly supported by Marie Curie IRG Grant and ERC Start Grant.}

\begin{document}

\begin{abstract}
I will talk about my recent work with Fernando Marques where we used  Almgren--Pitts Min-max Theory to settle some open questions in Geometry:  The Willmore conjecture, the Freedman--He--Wang conjecture for links (jointly with Ian Agol), and the existence of infinitely many minimal hypersurfaces in manifolds of positive Ricci curvature. Some open questions are suggested in the last section.
\end{abstract}

\maketitle

\section{Introduction}
I will start by introducing three problems in Geometry which, while quite distinct, have in common the fact that their solution comes from understanding unstable critical points in the space of all embedded hypersurfaces on a given manifold.  A panoramic  overview of variational methods in  Geometry can be found in the contribution of Fernando C. Marques \cite{marques-icm}.

\subsection{Willmore Conjecture}\label{willmore.section}

A central question in Mathematics has been the search for the ``optimal'' representative within a certain class of objects. Partially motivated by this principle, Thomas Willmore started in the $60$'s the quest for the ``optimal''  immersion of a surface in space. 

With that in mind, he associated to every compact surface $\Sigma\subset \R^3$ the quantity (now known as the {\em Willmore energy}),
\[\mathcal W(\Sigma)=\int_{\Sigma}\left(\frac{k_1+k_2}{2}\right)^2d\mu,\]
where $k_1,k_2$ are the principal curvatures of $\Sigma$. 

The Willmore energy is invariant under rigid motions, scaling,  and is large when the surface contains  long thin tubes or long thin holes thus detecting how ``bended'' the surface $\Sigma$ is in space. Less obvious, is the fact that the Willmore energy is also invariant under the inversion $x\mapsto x/|x|^{2}$ and thus invariant under conformal transformations. Willmore himself only found this  some years later  but, as  we explain soon, this was know already since the twenties.

It is worthwhile to remark that in applied sciences the Willmore energy had already made its appearance a long time ago,  under the name of {\em bending energy}, in order to study vibrating properties of thin plates.  In $1810$'s, Marie-Sophie Germain proposed, as the bending energy of a thin plate, the integral  with respect to the surface area of an even, symmetric function of the principal curvatures, in which case the Willmore energy is the simplest possible example (excluding the area). Similar quantities were also considered by Poisson around the same time.  

Moreover, the Willmore energy had also appeared in Mathematics through the work of Thomsen \cite{thomsen} and Blaschke \cite{blaschke} in the $1920$'s but their findings were forgotten and only brought to light after the interest in the Willmore energy increased.  In particular, Thomsen and Blaschke  were already aware of the conformal invariance of the Willmore energy.

Back to Willmore's quest for the ``best'' possible immersion, he  showed that round spheres have the least possible Willmore energy among all compact surfaces in space.  More precisely,  ever compact surface $\Sigma\subset \R^3$ has 
$$\mathcal W(\Sigma)\geq 4\pi$$ with equality only for  round spheres.

Having found the compact surface with least possible energy, he tried to find the torus in space with smaller energy than any other tori. It is interesting to note that, just by looking at the shape of tori in space, no obvious candidate stands out. Hence, Willmore  fixed a circle on a plane and considered tubes $\Sigma_r$ of a constant radius $r$ around  that circle. When $r$ is very small, $\Sigma_r$ is a  thin tube around the planar circle and thus its energy $\mathcal W(\Sigma_r)$ will be very large. If we keep increasing the value of $r$, the ``hole'' centered at the axis of revolution of the torus decreases and  eventually disappears for some $r_0$. Thus $\mathcal W(\Sigma_r)$ will be arbitrarily large for $r$ close  to $r_0$. Therefore $\mathcal W(\Sigma_r)$  must have an absolute minimum as $r$ ranges from $0$ to $r_0$, which Willmore computed to be $2\pi^2$. 

Up to scaling, the ``optimal'' torus that Willmore found has generating circle with radius 1 and center at distance $\sqrt{2}$ from  the axis of revolution:
$$
(u,v) \mapsto \big( (\sqrt{2} +\cos\, u) \cos\, v, (\sqrt{2}+\cos\,u)\sin\,v, \sin\,u) \in \mathbb{R}^3.
$$
In light of his findings, Willmore conjectured \cite{willmore}:
\begin{Wconj} 
Every compact surface $\Sigma$ of genus one has
$$\int_{\Sigma}\left(\frac{k_1+k_2}{2}\right)^2d\mu\geq 2\pi^2.$$
\end{Wconj}

It seems at first rather daring to  make such  conjecture after having tested it  only on a very particular one parameter family of tori. On the other hand, the torus Willmore found is special and had already appeared in Geometry: Inside the unit $3$-sphere $S^3$ in $\R^4$ there is a highly symmetric torus, the Clifford torus, which is given by $S^1(\frac{1}{\sqrt 2})\times S^1(\frac{1}{\sqrt 2})$. There is a stereographic projection  from the $3$-sphere onto space that sends the Clifford torus to the  ``optimal''  torus found by Willmore.

The richness of the  Willmore conjecture derives partially from the fact that the Willmore energy is invariant under conformal maps. One immediate consequence is that the conjecture can be restated for surfaces in the unit $3$-sphere $S^3$. Indeed, if $\Sigma$ is a compact surface in $S^3$ and $ \tilde \Sigma$ its image in $\R^3$ under stereographic projection, then one has
$$ \mathcal W(\tilde \Sigma)=\int_{\Sigma}1+ \left(\frac{ k_1+ k_2}{2}\right)^2d\mu,$$
where $ k_1,  k_2$ are the principal curvatures of $\Sigma$ {\em with respect} to the standard metric on $S^3$. For this reason, one  calls  the left-hand side of the above equation  the {\em Willmore energy}  $\mathcal W(\Sigma)$ of $\Sigma \subset S^3$. 

The conjecture had been verified in many special cases borrowing inspiration  from several distinct areas such as integral geometry, algebraic geometry, minimal surfaces, analysis or conformal geometry. We refer the reader to \cite{marques-neves} for the  history of partial results and mention simply the ones relevant to our work. 

In $1982$ it was proven by Li and Yau \cite{li-yau}  that the Willmore energy of any non-embedded surface must be at least $8\pi$ (which is strictly bigger than $2\pi^2$. In particular, it suffices to check the Willmore conjecture for embedded tori.  

Ros \cite{ros} proved in $1999$ the Willmore conjecture for tori in $S^3$ that are preserved by the antipodal map and his method motivated our approach.

Curiously, two biologists, Bensimon and Mutz \cite{bensimon}, verified the Willmore conjecture with the aide of a microscope while studying the physics of membranes!  They produced toroidal vesicles in a laboratory and observed that they assumed the  shape, which according to the Helfrich model \cite{helfrich} should be the minimizer for the Willmore energy,  of  the Clifford torus or its conformal images.

Jointly with Fernando Marques \cite{marques-neves} we showed that
\begin{thm}\label{willmore}
Every embedded compact surface $\Sigma$ of $S^3$ with positive genus has
$$\mathcal W(\Sigma)\geq 2\pi^2.$$
Equality only holds, up to rigid motions,  for  stereographic projections  of the Clifford torus.
\end{thm}

The rigidity statement characterizing the equality case in Theorem \ref{willmore} is optimal because, as we have mentioned, the Willmore energy is conformal invariant.

Using the  Li--Yau result  previously mentioned we obtain
\begin{coro}
 The Willmore conjecture holds.
 \end{coro}

\subsection{Energy of links}
The second application comes from the theory of links in $\mathbb{R}^3$. Let $\gamma_i: S^1 \rightarrow \R^3$, $i=1,2$, be a 2-component link, i.e., a pair of  closed curves in Euclidean three-space  with $\gamma_1(S^1) \cap \gamma_2(S^1) = \emptyset$. 

A $2$-component link is said to be {\em nontrivial} if it cannot be deformed without intersecting itself into two curves contained in disjoint balls. To every link $(\gamma_1,\gamma_2)$  one associates an  integer invariant, called the linking number ${\rm lk}(\gamma_1,\gamma_2)$, that intuitively measures how many times each curve winds around the other.

To every  $2$-component link $(\gamma_1,\gamma_2)$, O'Hara \cite{o'hara} associated an energy, called the  {\it M\"{o}bius cross energy}. Its definition is reminiscent of the  electrostatic potential energy and is given by
(\cite{o'hara}, \cite{freedman-he-wang}):
$$
E(\gamma_1,\gamma_2) = \int_{S^1 \times S^1} \frac{|\gamma_1'(s)||\gamma_2'(t)|}{|\gamma_1(s)-\gamma_2(t)|^2}\, ds\, dt.
$$ 
Freedman, He, and Wang studied this energy in detail and found that it has  the remarkable property of being invariant under conformal transformations of $\R^3$ \cite{freedman-he-wang}, just like the Willmore energy.  

Using Gauss formula for the linking number, one can see that  $$E(\gamma_1,\gamma_2)\geq 4\pi |{\rm lk}(\gamma_1,\gamma_2)|$$ and so  it is then natural to search for optimal configurations, i.e., minimizers of the M\"{o}bius energy.  
 This question can be given the following nice physical interpretation (see \cite{o'hara}). Assuming that each curve in the link  is non-conductive, charged uniformly and subject to a Coulomb's repulsive force, the equilibrium configuration the link will assume  should minimize the M\"{o}bius energy.
 
Freedman, He and Wang \cite{freedman-he-wang}  considered this question and after  looking at the particular the case where one of the link components is a planar circle, they made the following conjecture.

\begin{lconj}  The M\"{o}bius energy is minimized, among the class of all nontrivial links in $\R^3$, by the stereographic projection of the standard Hopf link in $S^3$.
\end{lconj}
 
The standard Hopf link $(\hat{\gamma}_1,\hat{\gamma}_2)$ in $S^3$  is described by
$$\hat{\gamma}_1(s)=(\cos s, \sin s,0,0) \in S^3\quad\mbox{and}\quad\hat{\gamma}_2(t)=(0,0,\cos t,\sin t) \in S^3,$$ and it is simple to check that  $E(\hat{\gamma}_1,\hat{\gamma}_2)=2\pi^2$.

In a joint work with Ian Agol and Fernando Marques \cite{agol-marques-neves} we showed that:

\begin{thm}\label{links} Let $(\gamma_1,\gamma_2)$ be a 2-component link in $\R^3$ with $|{\rm lk}(\gamma_1,\gamma_2)| = 1$.  Then $E(\gamma_1,\gamma_2) \geq 2\pi^2$. 

Equality only holds, up to rigid motions and orientation,  for  stereographic projections  of the Hopf link.
\end{thm}

It follows from a result of He \cite{He02} that it suffices to prove the conjecture for links $(\gamma_1,\gamma_2)$ with linking number ${\rm lk}(\gamma_1,\gamma_2)=\pm 1$.  Thus, we obtained the following corollary

\begin{coro} The  conjecture made by Freedman, He, and Wang holds.
\end{coro}

\subsection{Existence of embedded minimal hypersurfaces}

A question lying at the core of Differential Geometry,  asked  Poincar\'e \cite{poincare} in $1905$, is whether every closed Riemann surface always admits a closed geodesic.  

If the surface is not simply connected then we can minimize length in a nontrivial homotopy class and produce a closed geodesic. Therefore the question becomes considerably more interesting on a two-sphere, and the first breakthrough was in $1917$, due to Birkhoff \cite{birkhoff},   who found a closed geodesic for any  metric on a two-sphere.  

Later, in a remarkable work, Lusternik and Schnirelmann \cite{lusternik}  showed that every metric on a $2$-sphere admits three simple (embedded) closed geodesics (see also \cite{ballman, grayson, jost2, klingenberg, lusternik2, taimanov}).  This result is optimal because there are ellipsoids which admit no more than three simple closed geodesics.

This suggests the question of whether we can find  an infinite number of geometrically distinct closed geodesics in any closed surface. It is not hard to find infinitely many closed geodesics
when the genus of the surface is positive.

The case of the sphere was finally settled in $1992$ by Franks \cite{franks} and Bangert \cite{bangert}. Their  works combined  imply that every metric on a two-sphere admits an infinite number of closed geodesics. Later, Hingston \cite{hingston} estimated the number of closed geodesics of length at most $L$ when $L$ is very large.

Likewise, one can ask whether every closed Riemannian manifold admits a closed minimal hypersurface. When the ambient manifold has topology one can find minimal hypersurfaces by minimization and so, like in the surface case, the question is more challenging when every hypersurface is homologically trivial. Using min-max methods, and building on earlier work of Almgren, Pitts \cite{pitts} in $1981$ proved that  every compact Riemannian $(n+1)$-manifold with $n\leq 5$ contains a smooth, closed, embedded  minimal hypersurface. One year later, Schoen and Simon \cite{schoen-simon} extended this result to any dimension, proving the existence of a closed,  embedded minimal hypersurface with a singular set of Hausdorff codimension at least $7$.

When $M$ is diffeomorphic to a $3$-sphere, Simon--Smith \cite{smyth} showed the existence of  a minimal embedded sphere using min-max methods (see also \cite{colding-delellis}).

Motivated by these results, Yau made the following conjecture \cite{yau1} (first problem in the Minimal Surfaces section):
\begin{yconj} Every compact $3$-manifold $(M,g)$ admits an infinite number of  smooth, closed, immersed minimal surfaces.
\end{yconj}

Lawson \cite{lawson70} showed in $1970$ that  the round $3$-sphere admits embedded minimal surfaces of every possible genus. 

When $M$ is a compact hyperbolic $3$-manifold, Khan and Markovic \cite{kahn-markovic} found an infinite number of incompressible surfaces in $M$ of arbitrarily high genus. One can then minimize energy in their homotopy class and obtain an  infinite number of  smooth, closed, immersed minimal surfaces.

Jointly with Fernando Marques \cite{marques-neves-3} we showed
\begin{thm}\label{ythm} Let $(M^{n+1},g)$ be  a compact Riemannian manifold with $2\leq n\leq 6$ and a metric of positive Ricci curvature.

Then $M$ contains an infinite number of distinct, smooth, embedded, minimal hypersurfaces.
\end{thm}

Until Theorem \ref{ythm} was proven, it was not even known whether metrics on the $3$-sphere arbitrarily close to the round metric also admit an infinite number of minimal  surfaces.  

I find a fascinating problem to shed some light into the asymptotic behaviour of the minimal surfaces given by Theorem \ref{ythm}.

\section{Almgren--Pitts Min-max Theory}

As  mentioned in the Introduction, Theorem \ref{willmore}, Theorem \ref{links}, and Theorem \ref{ythm}, follow from understanding   the topology of the space of all embedded hypersurfaces.

In very general terms, the guiding principle of Morse Theory is that given a space and a function defined on that space, the topology of the space forces the function to have certain critical points. For instance, if the space has a  $k$-dimensional nontrivial cycle, then the function should have   a critical point of index at most $k$.

The space we are interested  is $\mathcal Z_n(M)$,  the space of all orientable compact hypersurfaces with possible multiplicities in a compact Riemannian $(n+1)$-manifold $(M,g)$ with $n\geq 2$. If we allow for non-orientable hypersurfaces as well, the space is denoted by $\mathcal Z_n(M;\Z_2)$. 

 These spaces are studied in the context of Geometric Measure Theory and come  with a well understood topology (flat topology) and  equipped with a natural functional  which associates to every element in  $\mathcal Z_n(M)$ (or $\mathcal Z_n(M;\Z_2)$) its $n$-dimensional volume.  

I  will try to keep the discussion with as little technical jargon as possible in order to convey the main ideas and thus ignore almost all technical issues.  

Critical points of the volume functional are called {\em minimal hypersurfaces} and their index is the number of independent deformations that decrease the area. For instance, on the $3$-torus with the flat metric, there is a natural flat $2$-torus which minimizes area in its homology class and thus it is a minimal surface of index zero. Likewise, on the $3$-sphere with the round metric, the equator (which has area $4\pi$) is a minimal sphere and has the property that if we ``push'' it up into the northern hemisphere then its area is decreased. Hence, its index is at least one. On the other hand, one can check that a deformation of the equator that preserves the enclosed volume is never area decreasing. Thus the equator has index one.

Almgren \cite{almgren-varifolds} started in the $60$'s the study of Morse Theory for the volume functional  on $\mathcal Z_n(M)$ (or $\mathcal Z_n(M;\Z_2)$) and that continued through the $70$'s  jointly with Pitts, his Phd student.   I present now the basic principles of Almgren--Pitts Min-max Theory.

Suppose $X$ is a topological space and $\Phi:X\rightarrow \mathcal Z_n(M)$ a continuous function. Consider
$$[\Phi]=\{\Psi:X\rightarrow \mathcal Z_n(M):\, \Psi\mbox{ homotopic to }\Phi\mbox{ relative to }\partial X\}.$$
Note that if $\Psi\in [\Phi]$ then $\Phi=\Psi$ on $\partial X$.
To the homotopy class $[\Phi]$ we associate the number, called the {\em width},
$${\bf L}([\Phi])=\inf_{\Psi\in[\Phi]}\sup_{x\in X}{\vol}(\Psi(x)).$$

The Almgren-Pitts Min-max Theorem \cite{pitts} can be stated as
\begin{thm}[Min-max Theorem] \label{min.max.thm.light}Assume that ${\bf L}([\Phi])>\sup_{x\in \partial X}{\vol}(\Phi(x))$.

There is a compact embedded minimal hypersurface $\Sigma$ (with possible multiplicities) such that 
\begin{equation*}\label{almgren.pitts.formula}
{\bf L}([\Phi])={\vol }(\Sigma).
\end{equation*}
\end{thm}

The theorem also holds for $\mathcal Z_n(M;\Z_2)$ with no modifications.

The support of $\Sigma$ is smooth outside a set of codimension $7$ and thus smooth if $n\geq 6$ (for $n\geq 5$ the regularity theory was done by Schoen and Simon \cite{schoen-simon}). The Min-max Theorem allows for $\Sigma$  to be a union of disjoint hypersurfaces, each  with some multiplicity. More precisely 
$$\Sigma=n_1\Sigma_1+\ldots+n_k\Sigma_k,$$ where $n_i\in \N$, $i=1,\ldots,k$, and $\{\Sigma_1,\ldots, \Sigma_k\}$ are embedded minimal surfaces with disjoint supports. 

Naturally, if the space of parameters $X$ is a $k$-dimensional, we expect the index of $\Sigma$ to be at most $k$ but  this fact has not been proven.

The condition ${\bf L}([\Phi])>\sup_{x\in \partial X}{\vol}(\Phi(x))$  means  that $[\Phi]$ is capturing some nontrivial topology of  $\mathcal Z_n(M)$. The guiding philosophy behind Min-max Theory consists in  finding examples of homotopy classes satisfying this condition and then use the Min-max Theorem  to deduce geometric consequences.

The next example illustrates well this methodology. Given $f:M\rightarrow [0,1]$ a Morse function,  consider the continuous map 
$$\Phi:[0,1]\rightarrow \mathcal Z_n(M),\quad \Phi(t)=\partial\{x\in M:f(x)<t\}.$$
We have $\Phi(0)=\Phi(1)=0$ because all elements in $\mathcal Z_n(M)$ with zero volume are identified to be the same and Almgren showed in \cite{almgren} that $L([\Phi])>0$, i.e., $[\Phi]$ is a nontrivial element of $\pi_1(\mathcal Z_n(M),\{0\})$. Using Min-max Theorem one  obtains the existence of a minimal hypersurface in $(M,g)$ whose volume realizes  $L([\Phi])$ and thus $(M,g)$ admits a minimal embedded hypersurface which is smooth outside a set of codimension $7$. This application was one of the motivations for Almgren and Pitts to develop their Min-max Theory.  

\section{The $2\pi^2$ Theorem}

Let $I^k$ denote a closed $k$-dimensional cube and $B_r(p)$ denote the geodesic ball in $S^3$ of radius $r$ centered at $p$. 

We present  a criteria due to  Marques and myself  \cite{marques-neves}  to ensures that a map $\Phi:I^5\rightarrow \mathcal Z_2(S^3)$ determines a nontrivial $5$-dimensional homotopy class in $\mathcal{Z_2(S^3)}.$ 

Let $\mathcal G_o$ be the set of all oriented geodesic spheres in $\mathcal Z_2(S^3)$. Each nonzero element  in $\mathcal G_o$ is determined by its center and radius. Thus this space is homeomorphic to $S^3\times[-\pi,\pi]$ with an equivalence relation that identifies $S^3\times \{-\pi\}$ and $S^3\times \{\pi\}$ all with the zero in $\mathcal Z_2(S^3)$. 

The maps $\Phi$ we consider have the property that 
$$\Phi(I^4\times\{1\})=\Phi(I^4\times\{0\})=\{0\}\quad\mbox{and}\quad \Phi(I^4\times I)\subset \mathcal G_o.$$
Hence $\Phi(I^5)$ can be thought of a $5$-cycle in $\mathcal Z_2(S^3)$ whose boundary lies in $\mathcal G_o$ and thus $[\Phi]$ can be seen as an element of $\pi_5(\mathcal Z_2(S^3), \mathcal G_o)$. The next theorem gives a condition under which $[\Phi]\neq 0$ in  $\pi_5(\mathcal Z_2(S^3), \mathcal G_o)$, i.e., the image of $\Phi$ cannot be homotoped into the set of all  geodesic spheres.

\begin{thm}\label{thmd} Let $\Phi:I^5 \rightarrow 
 \mathcal Z_2(S^3)$ be a continuous map  such that:
\begin{itemize}
\item[(1)] $\Phi(x,0)=\Phi(x,1) = 0$ for any $x\in I^4$;
\item[(2)] for any $x\in \partial I^4$ fixed we can find $Q(x)\in S^3$ such that
$$\Phi(x,t)=\partial B_{{\pi}t}(Q(x)),\quad 0\leq t\leq 1.$$
In particular, $\Phi(I^5)\subset \mathcal G_o$.
\item[(3)] the center map $Q: \partial I^4\rightarrow S^3$ has ${\rm deg}(Q) \neq 0$.
\end{itemize}
Then 
$${\bf L}([\Phi])>4\pi=\sup_{x\in \partial I^5} {\bf M}(\Phi(x)).$$
\end{thm}

Condition (3) is crucial, as the next example shows. Consider
$$ \Phi:I^5\rightarrow \mathcal{Z}_2(S^3), \quad \Phi(x,t)=\partial B_{\pi t}(p),$$
where $p$ is a fixed point in $S^3$. Conditions (1) and (2)  of Theorem \ref{thmd} are satisfied but  ${\bf L}([\Phi])=4\pi$ because $\Phi(I^5)\subset \mathcal{G}_o$.

\begin{proof}[Sketch of proof]
The idea  for the proof is, in very general terms, the following.  Let $\mathcal R\subset \mathcal G_o$ denote the space of all oriented great spheres (which is homeomorphic to $S^3$). With this notation, note that from condition (2) we have that $$\Phi(\partial I^5\times\{1/2\})\subset \mathcal R.$$
Moreover, the map $\Phi:\partial I^5\times\{1/2\}\rightarrow \mathcal R\approx S^3$ has degree equal to  ${\rm deg}(Q)$ and thus nonzero by condition (3). For simplicity, suppose we can find $\Psi\in[\Phi]$ so that
 $${\bf L}([\Phi])=4\pi=\sup_{x\in I^5}{\rm area}(\Psi(x)).$$ 
In particular, for any given continuous path $\gamma: [0, 1]\rightarrow I^5$ connecting $I^4 \times \{0\}$ to $I^5 \times \{1\}$, $\Psi\circ\gamma$ is optimal as a one-parameter sweepout of $S^3$ and thus  it must intersect $\mathcal R$, i.e., contain a great sphere.
 
 Then  $R=\Psi^{-1}(\mathcal R)$ should be  a $4$-dimensional cycle in $I^5$ separating bottom from top and with $\partial R=\partial I^4\times \{1/2\}$. Hence
$$\Psi_{*}[\partial R]=\partial[ \Psi(S)]=\partial[\mathcal R]=0\mbox{ in }H_3(\mathcal R,\Z).$$
On the other hand, $\Phi=\Psi$ on $\partial S=\partial I^4\times \{1/2\}$ and so
$$\Psi_{*}[\partial R]=\Phi_{*}[\partial I^4\times\{1/2\}]={\rm deg}(Q)[\mathcal R]\neq 0$$
and this is a contradiction.
\end{proof}
Suppose now that $\Phi$ is a map satisfying the hypothesis of Theorem \ref{thmd}. From the Min-max Theorem we obtain the existence of $\Sigma$, an embedded minimal surface, such that ${\bf L}([\Phi])={\rm area}(\Sigma)>4\pi$. Moreover, it is natural to expect that $\Sigma$ has index  at most $5$ because we are dealing with a $5$-parameter family of surfaces.   

Urbano \cite{urbano} in $1990$ classified minimal surface of $S^3$ with low index and he gave a rather elegant and short proof of

\begin{thm}[Urbano's Theorem]
\label{urbano1} Assume  $S$ is a closed  embedded  minimal surface in $S^3$ having  ${\rm index}(S)\leq 5.$ 

Then, up to ambient isometries, $S$ is either a great sphere (index one) or the Clifford torus (index five) up to ambient isometries.
\end{thm}

\begin{remark}
We already argued that the great sphere has index one. The Clifford torus has index five because  unit speed normal deformations decreases area, the four parameter space of conformal dilations (to be seen later) also decrease area, and these five deformations are linearly independent. 
\end{remark}

Going back to our discussion, we see that $\Sigma$ cannot be a great sphere because its area is ${\bf L}([\Phi])>4\pi$ and so it has to be a Clifford torus with area $2\pi^2$. This heuristic discussion motivates the
\begin{thm}[$2\pi^2$ Theorem] Assume that $\Phi$ satisfied the hypothesis of Theorem \ref{thmd}. Then
$$\sup_{x\in I^5}{\rm area}(\Phi(x))\geq 2\pi^2.$$
\end{thm}
A  proof can be found in \cite{marques-neves}. Because the Almgren--Pitts theory does not provide us with the fact that the index of $\Sigma$ is at most $5$, we had to use a new set of arguments to prove the index estimate in the case we were interested. 

\section{Strategy to prove Theorem \ref{willmore}}

We sketch the proof of the inequality in  Theorem \ref{willmore}. The complete argument can be found in \cite{marques-neves}. 

The conformal maps of $S^3$ (modulo isometries) can be parametrized by the open unit $4$-ball $B^4$, where to each nonzero $v\in B^4$ we consider the conformal dilation $F_v$ centered at $\frac {v}{|v|}$ and $-\frac{v}{|v|}$. Composing with the stereographic  projection $\pi:S^3\setminus \{-\frac {v}{|v|}\} \rightarrow \mathbb{R}^3$, the map $$\pi\circ F_v\circ \pi^{-1}:\R^3\rightarrow \R^3$$ corresponds to a dilation in space centered at the origin, where the dilation factor tends to infinity as $|v|$ tends to one.

Given a compact embedded surface $S\subset S^3$ and $-\pi\leq t\leq\pi$, we denote by $S_t$ the surface at distance $|t|$ from $S$, where $S_t$ lies in the exterior (interior) of $S$ if $ t \geq 0$ ($t\leq 0$). Naturally, $S_t$ might not be a smooth embedded surface due to the existence of possible focal points but it will always be well defined in the context of Geometric Measure Theory.

We can now define the following $5$-parameter family $\{\Sigma_{(v,t)}\}_{(v,t)\in B^4\times [\pi,\pi]}$ of surfaces in $S^3$ given by
$$ \Sigma_{(v,t)}=(F_v(\Sigma))_t\in \mathcal Z_2(S^3).$$
One crucial property of this $5$-parameter family is the following.
\begin{thm}[Heintze--Karcher Inequality] For every $(v,t)\in  B^4\times [\pi,\pi]$ we have
$${\rm area}(\Sigma_{(v,t)})\leq \mathcal W(\Sigma).$$
\end{thm}
A related result was proven by Ros in \cite{ros}.

In order to apply the $2\pi^2$ Theorem it is important that we understand the behaviour $\Sigma_{(v,t)}$ as $(v,t)$ approaches the boundary of $B^4\times [\pi,\pi]$. The fact that the diameter of $S^3$ is $\pi$ implies that $\Sigma_{(v,\pm\pi)}=0$ for all $v\in B^4$ and so we are left to analyze what happens when $v$ approaches $S^3$.

Assume $v$ in the $4$-ball tends to $p\in S^3$. If $p$ does not belong to $\Sigma$, then  $F_v(\Sigma)$  is ``pushed'' into  $\{-p\}$ as $v$ tends to $p$ and so ${\rm area}(F_v(\Sigma))$ tends to zero. When $p$ lies in $\Sigma$ the situation is considerably more subtle. Indeed,  if $v$ approaches $p$ radially, i.e., $v=sp$ with $0<s<1$, then $F_{sp}(\Sigma)$ converges, as $s$ tends to $1$, to the unique great sphere tangent to $\Sigma$ at $p$. Thus the continuous function in $S^3$ given by $p\mapsto {\rm area}\,(\Sigma_{sp})$ tends, as $s\to 1$, to a discontinuous function that is zero outside $\Sigma$ and $4\pi$ along $\Sigma$. Hence, for any $0<\alpha<4\pi$, there must exist a sequence $\{v_i\}_{i\in\N}$ in  $B^4$ tending to $\Sigma$ so  that $ {\rm area}\,(\Sigma_{v_i})$ tends to $\alpha$ and so it is natural to expect that the convergence of $F_v(\Sigma)$ depends on how $v$ approaches $p\in \Sigma$. A careful analysis revealed that, depending on the angle at which $v$ tends to $p$, $F_v(\Sigma)$ tends to a geodesic sphere tangent to $\Sigma$ at $p$, with radius and center depending on the angle of convergence.

Initially this behaviour was a source of perplexity but then we realized that, even if the parametrization was becoming discontinuous near the boundary of the parameter space, the closure of the family $\{\Sigma_{(v,t)}\}_{(v,t)\in B^4\times [\pi,\pi]}$ in $\mathcal Z_2(S^3)$ was a ``nice'' continuous $5$-cycle in $\mathcal Z_2(S^3)$.   Hence, we were able to  reparametrize this family and obtain a  continuous map $\Phi:I^5\to \mathcal Z_2(S^3)$ with $\Phi(I^5)$  equal to the closure of $\{\Sigma_{(v,t)}\}_{(v,t)\in B^4\times [\pi,\pi]}$ in $\mathcal Z_2(S^3)$ and satisfying conditions (1) and (2) of Theorem \ref{thmd}. 

Finally, and most important of all, we showed that the degree of the center map $Q$ in condition (2) of Theorem \ref{thmd} is exactly the genus of $\Sigma$. This point is absolutely crucial because it showed us that the map $\Phi$ ``remembers'' the genus of the surface $\Sigma$. Thus, when the genus  is positive, condition (3) of Theorem \ref{thmd} is also satisfied and we obtain from the $2\pi^2$ Theorem that 
$$2\pi^2\leq  \sup_{x\in I^5}{\rm area}((\Phi(x)).$$
On the other hand, because $\Phi(I^5)$ is equal to the closure of $\{\Sigma_{(v,t)}\}_{(v,t)\in B^4\times [\pi,\pi]}$ in $\mathcal Z_2(S^3)$, we obtain from the Heintze--Karcher Inequality that 
$$\sup_{x\in I^5}{\rm area}(\Phi(x))\leq \mathcal W(\Sigma).$$
This means that $\mathcal W(\Sigma)\geq 2\pi^2$, which is the statement we wanted to prove. 

\section{Strategy to prove Theorem \ref{links}}

The approach to prove Theorem  \ref{links} is similar to the one used in Theorem \ref{willmore}.   The conformal invariance of  the energy implies that  it suffices to consider links $(\gamma_1,\gamma_2)$ in $S^3$. For each link $(\gamma_1,\gamma_2)$ in $S^3$  we construct 
a suitable family
$\Phi:I^5 \rightarrow \mathcal{Z}_2(S^3)$  that satisfies conditions (1) and (2) of Theorem \ref{thmd}. Moreover, we will also show that if $|{\rm lk}(\gamma_1,\gamma_2)| = 1$ then condition (3) of Theorem \ref{thmd} is will also be satisfied.  Hence we can apply the $2\pi^2$ Theorem  and conclude that 
$$\sup_{x\in I^5}{\rm area}(\Phi(x))\geq 2\pi^2.$$ 
On the other hand, the map $\Phi$  is constructed
so that ${\rm area}(\Phi(x)) \leq E(\gamma_1,\gamma_2)$ for each $x \in I^5$ and this implies the inequality in Theorem \ref{links}.

We give a brief indication of how the map $\Phi$ is constructed.  

To every pair of curves in $\R^4$ there is a natural way to construct a ``torus'' in $S^3$. More precisely,
 given two curves $(\gamma_1,\gamma_2)$ in $\R^4$, the {\em Gauss map}  is denoted by
$$G(\gamma_1,\gamma_2):S^1 \times S^1 \rightarrow S^3, \quad (s,t)\mapsto \frac{\gamma_1(s)-\gamma_2(t)}{|\gamma_1(s)-\gamma_2(t)|}
$$
and we consider $G(\gamma_1,\gamma_2)_{\#}(S^1\times S^1)$ in $\mathcal Z_2(S^3)$.  Furthermore, one can check that
\begin{equation*}\label{mass.ineq.intro}
{\rm area}(G(\gamma_1,\gamma_2)_{\#}(S^1\times S^1))\leq E(\gamma_1,\gamma_2).
\end{equation*}
For instance, if $(\gamma_1,\gamma_2)$ is the Hopf link then $G(\gamma_1,\gamma_2)_{\#}(S^1\times S^1)$ is the Clifford torus and the inequality above becomes an equality.
 
Given $v\in B^4$, we consider the conformal map $F_v$ of $\R^4$ given by an inversion centered at $v$. The conformal map $F_v $ sends the unit $4$-ball $B^4$ to some other ball centered at $c(v)=\frac{v}{1-|v|^2}$. We consider $$g:B^4\times (0,+\infty)\rightarrow \mathcal Z_2(S^3)$$ given by 
$$g(v,z)=G\left(F_v\circ\gamma_1,\lambda(F_{v}\circ\gamma_2-c(v))+c(v)\right)_{\#}(S^1\times S^1).$$
Intuitively, $g(v,z)$ is the image of  the Gauss map of the link obtained by applying the conformal transformation $F_v$ to  $(\gamma_1,\gamma_2)$ and then dilating the  curve $F_v\circ\gamma_2$ with respect to the center $c(v)$ by a factor of $\lambda$. Note that both curves $F_v\circ\gamma_1$ and $\lambda(F_{v}\circ\gamma_2-c(v))+c(v)$ are contained in spheres centered at $c(v)$.  

The $5$-parameter family we just described also enjoys a Heintze--Karcher type-inequality, meaning that for all $(v,\lambda)\in B^4\times (0,+\infty)$ we have
$${\rm area}(g(v,z))\leq E(\gamma_1,\gamma_2).$$
The map $\Phi$ is constructed via a reparametrization of $g$.

\section{Gromov--Guth families} In order to apply the Min-max Theorem on a general manifold $M$, it is important that one understands the  homotopy groups of the space  $\mathcal Z_{n}(M;\Z_2)$. This was done by Almgren \cite{almgren} in $1962$ and  he showed that
$$\pi_1(\mathcal Z_{n}(M;\Z_2)) = \Z_2\quad\mbox{and}\quad \pi_i(\mathcal Z_{n}(M;\Z_2))=0\quad\mbox{ if } i>1.$$
 Thus $\mathcal Z_{n}(M;\Z_2)$ is weakly homotopic equivalent to  $\RP^{\infty}$ and so we should expect that   $\mathcal Z_{n}(M,\Z_2)$ contains, for every $p\in \N$, an homotopically nontrivial $p$-dimensional projective space.

 From the weak homotopy equivalence with $\RP^\infty$, we have that $$H^k(\mathcal Z_{n}(M;\Z_2))=\Z_2\quad\mbox{ for all $k\in\N$ with generator $\bar\lambda^k$.}$$ 
 We are interested in studying maps $\Phi:X\rightarrow \mathcal Z_{n}(M;\Z_2)$ whose image detects $\bar \lambda^p$ for some $p\in\N.$
 
Given a simplicial complex $X$, a continuous map $\Phi:X\rightarrow \mathcal Z_{n}(M;\Z_2)$ is called a {\em p-sweepout} if  $\Phi^*(\bar\lambda^p)\neq 0$ in $H^p(X;\Z_2)$. Heuristically, a continuous map $\Phi:X\rightarrow \mathcal Z_{n}(M;\Z_2)$ is called a {\em p-sweepout} if for every set $\{x_1,\ldots,x_p\}\subset M$, there is $\theta\in X$ so that $\{x_1,\ldots,x_p\}\subset \Phi(\theta)$.

Gromov \cite{gromov0,gromov,gromov2} and Guth \cite{guth} studied $p$-sweepouts of $M$. 

We now check that $p$-sweepouts exists for all $p\in M$. Let $f\in  C^{\infty}(M)$ be a Morse function and consider the  map 
$${\Phi}:\RP^p\rightarrow \mathcal{Z}_n(M;\Z_2),$$
given by
 $$
{\Phi}([a_0,\ldots,a_p])=\partial\left \{x\in M:a_0+a_1f(x)+\ldots+a_p f^p(x)<0\right\}.$$
Note that the map $\Phi$ is well defined  because opposite orientations on the same hypersurface determine the same element in $\mathcal{Z}_n(M;\Z_2).$  A typical element of ${\Phi}([a_0,\ldots,a_p])$ will  be $f^{-1}(r_1)\cup\ldots\cup f^{-1}(r_j)$ where $r_1,\ldots, r_j$ are the real roots of the polynomial $p(t)=a_0+a_1t+\ldots+a_pt^p$. It  easy to see that $\Phi$ satisfies the heuristic definition of a $p$-sweepout given above and in \cite{marques-neves-3} we check that the map is indeed a $p$-sweepout.

Denoting the set of  all $p$-sweepouts of $M$ by $\mathcal P_p$, the {\em p-width} of $M$ is defined 
{as $$\omega_p(M)=\inf_{\Phi \in \mathcal P_p}\sup\{{\vol }(\Phi(x)): x\in {\rm dmn}(\Phi)\},$$
where ${\rm dmn}(\Phi)$ is the domain of $\Phi$.} 

 It is interesting to compare the $p$-width with the min-max definition of the $p^{th}$-eigenvalue of $(M,g)$.  Set $V= W^{1,2}(M)\setminus\{0\}$  and recall that 
$$\lambda_p=\mathop{\inf}_{(p+1)-\mbox{plane }P\subset V}\max \left\{ \frac{\int_M|\nabla f|^2 dV_g}{\int_M f^2 dV_g}: f\in P\right\}.$$
Hence one can see $\{\omega_p(M)\}_{p\in\N}$ as a nonlinear analogue of the Laplace spectrum of $M$, as proposed by  Gromov \cite{gromov0}. 

The asymptotic behaviour of $w_p(M)$ is governed by the following result proven by Gromov \cite{gromov0} and Guth \cite{guth}.
\begin{thm}[Gromov and Guth's Theorem] There exists a positive constant $C=C(M,g)$ so that,  for every $p\in \mathbb{N}$, $$C^{-1}p^{\frac{1}{n+1}}\leq \omega_p(M) \leq C p^{\frac{1}{n+1}}.$$
\end{thm}
The idea to prove the lower bound is, roughly speaking, the following. Choose $p$ disjoint geodesic balls $B_1,\ldots, B_p$ with radius proportional to $p^{-\frac{1}{n+1}}$. For every $p$-sweepout $\Phi$ one can find  $\theta\in {\rm dmn}(\Phi)$ so that $\Phi(\theta)$ divides each geodesic ball into two pieces with almost identical volumes. Hence, when $p$ is sufficiently large,  the isoperimetric inequality implies that $\Phi(\theta)\cap B_i$ has volume no smaller than $c(n)p^{-\frac{n}{n+1}}$ for all $i=1,\ldots,p$, where $c(n)$ is a universal  constant.   As a result $\Phi(\theta)$ has volume greater or equal than $c(n)p^{\frac{1}{n+1}}$.

The upper bound  can be proven using a very nice bend-and-cancel argument introduced by Guth \cite{guth}. In the case when $M$ is a $n+1$-dimensional sphere  $S^{n+1}$, the upper bound has the following simple explanation. If we consider the set of all homogenous harmonic polynomials in $S^{n+1}$ with degree less or equal than $d\in \N$, we obtain a vector space of dimension $p(d)+1$, where $p(d)$ grows like $d^{n+1}$. Considering the zero set of all these polynomials we obtain a map $\Phi$ from a $p(d)$-dimensional projective plane  into $\mathcal Z_n(S^{n+1};\Z_2)$. Crofton formula implies that the zero-set of each of these polynomials  has volume at most $\omega_nd$, where $\omega_n$ is the volume of an $n$-sphere. Thus, for every  $\theta \in \RP^{p(d)}$, ${\vol}(\Phi(\theta))$ is at most a fixed multiple of $p(d)^{\frac{1}{n+1}}$.

\section{Strategy to prove Theorem \ref{ythm}}

The idea to find an infinite number of minimal surfaces consists in applying the Min-max Theorem to the family of $p$-sweepouts $\mathcal P_p$ for all $p\in\N$.

The first thing we show is that if $\omega_p(M)=\omega_{p+1}(M)$ for some $p\in \N$, then $M$ admits an infinite number of minimal embedded hypersurfaces. We achieve this using Lusternick-Schnirelman  and, roughly speaking, the idea is as follows (for details see \cite{marques-neves-3}): 

Suppose for simplicity that  $\omega_p(M)=\omega_{p+1}(M)=\sup_{x\in X}\vol(\Phi(x))$ for some $p+1$-sweepout $\Phi$. We argue by contradiction and assume that $\Omega$, the set of all embedded minimal hypersurfaces (with possible multiplicities) and volume at most $\omega_{p+1}(M)$, is finite. 

Let $K=\Psi^{-1}(\Omega)$ and $\lambda=\Phi^*(\bar\lambda)$, where $\bar \lambda$ generates $H^1(\mathcal Z_{n}(M;\Z_2);\Z_2)$. 

We must have $\lambda$ vanishing on $K$ because otherwise there would exist a curve $\gamma$ in $K$ so that $\lambda(\gamma)=\bar \lambda(\Phi\circ \gamma)\neq 0$, i.e. $\Phi\circ \gamma$ would be a  nontrivial element of $\pi_1(\mathcal Z_{n}(M;\Z_2))$. But $\Phi\circ \gamma$ has image contained in the finite set $\Omega$, which implies  it is constant, and thus contractible.

Therefore $\lambda^p$ cannot vanish on $X\setminus K$ because otherwise $\lambda^{p+1}=\lambda^p\smile \lambda$ would be zero on $(X\setminus K)\cup K=X$ and this is impossible because  $\Phi$ is a $p+1$-sweepout. As a result $\Phi_{|X\setminus K}$ is a $p$-sweepout whose image contains no minimal hypersurfaces (with possible multiplicities), and so we pull-tight the family to obtain another $p$-sweepout $\Psi$ so that 
$$
\omega_p(M)\leq\sup_{x\in X\setminus K}\vol(\Psi(x))<\sup_{x\in X\setminus K}\vol(\Phi(x))\leq \sup_{x\in X}\vol(\Phi(x))=\omega_p(M).$$ 
This gives us the desired contradiction.

Hence we can assume that the sequence $\{\omega_p(M)\}_{p\in\N}$ is strictly increasing. 

We then argue again by  contradiction and assume  that there exist only finitely many smooth, closed, embedded minimal hypersurfaces with multiplicity one, and we call this set $\Lambda$.  Using the Min-max Theorem we have that
$$\omega_p(M)=n_{p,1}\vol(\Sigma_{p,1})+\ldots+n_{p,k}\vol(\Sigma_{p,k}), \quad n_{p,1},\dots,n_{p,k}\in\N,$$
where $\{\Sigma_{p,1},\ldots,\Sigma_{p,k}\}$ are multiplicity one minimal hypersurfaces with disjoint support.
Because $M$ has positive Ricci curvature,   Frankel's Theorem \cite{frankel} says that any two minimal embedded hypersurfaces intersect  and so 
$\omega_p(M)=n_p\vol(\Sigma_p)$ for some $\Sigma_p\in \Lambda$ and $n_p\in\N$. The fact that  $\omega_p(M)$ is strictly increasing  and a counting argument shows that $\Lambda$ being finite implies that  $\omega_p(M)$ must grow linearly  in $p$.  This  is in contradiction with the sublinear growth of $\omega_p(M)$ in $p$ given by the Gromov and Guth's Theorem.

\section{Open problems}
Min-max Theory is an exciting technique which I think can be used not only to solve other open questions in Geometry but also to provide some new directions. Some of these questions are well-known and others arose from extensive discussions with Fernando Marques.

\subsection*{Min-max questions}
The Almgren-Pitts Min-max Theory  does not provide index estimates for the min-max minimal hypersurface. The importance of this issue was already clear to Almgren \cite{almgren-varifolds} who wrote\\

{\em ``The chief utility of the homology approach would lie in the attempt to assign a topological index to stationary integral varifolds in some analytically useful way.''}

 \medskip
 
A folklore conjecture states that the if the homotopy class in the Min-max Theorem is defined with $k$-parameters, then the minimal hypersurface given by the Min-max Theorem has index at most $k$. It  is implicit in the conjecture  that one finds a meaningful way of assigning an index to a minimal embedded hypersurface with multiplicities.
When $k=1$, the conjecture was confirmed by  Marques and myself if the ambient manifold is three dimensional and the metric has positive Ricci curvature \cite{marques-neves-3}. Later this was extended to the case where the ambient manifold has dimension between three and seven and the metric has positive Ricci curvature \cite{zhou}. 

Naively, one should also expect that for bumpy metrics the index is bounded from below by the number of parameters.

Many of the subtle issues in Min-max Theory are related with the fact that the min-max minimal hypersurface can have multiplicities. That said, one does not know an example where the width of some homotopy class is realized by an unstable minimal hypersurface with multiplicity. For instance,  does the equator with multiplicity two (and so area $8\pi$) realizes the width of some homotopy class in the round $3$-sphere? 
It is highly conceivable that unstable minimal surfaces with higher multiplicity can be approximated (in the varifold norm) by a sequence of embedded minimal surfaces with smaller area and this is one of the reasons the question is interesting.

\subsection*{Old questions} We now mention  four well known open problems which could be answered using Min-max Theory. 

The first two  are natural generalizations of the Willmore conjecture.  

The  Willmore conjecture in $S^4$ states  that among all tori in $S^4$, the Clifford torus minimizes the Willmore energy.  It is interesting that in this case there are minimal embedded projective planes (Veronese surface) which have area ($6\pi$) smaller than the Clifford torus and bigger than the equator.

 For higher genus surfaces, Kusner   \cite{kusner96} conjectured that the Lawson  minimal surface $\xi_{1,g}$ minimizes the Willmore energy among all surfaces of genus $g$ (numerical evidence was provided in \cite{hsu-kusner-sullivan}). It would be extremely interesting to find the index of $\xi_{1,2}$. Wishful thinking would  suggest $9$ but there is no real evidence.


The third problem consists of finding, among all  non-totally geodesic minimal hypersurfaces in the unit $n$-sphere $S^n$, the one with least possible volume. The conjecture, due to Solomon, is that these minimal hypersurface are given by
$$S^{m-1}\left(\sqrt\frac{m-1}{2m-1}\right)\times S^m\left(\sqrt\frac{m}{2m-1}\right)\subset S^{n}$$
 if $n=2m$, and by
 $$S^{m-1}\left(\frac{1}{\sqrt 2}\right)\times S^{m-1}\left(\frac{1}{\sqrt 2}\right)\subset S^{n}$$
if $n=2m-1$.  In $S^3$ this conjecture was confirmed by Marques and myself in \cite{marques-neves} and in the general case there was some progress due to White and Ilmanen \cite{white-ilmanen}. 

It would be desirable to have a sharp index characterization similar Urbano's Theorem  for  each of these three problems. In the third problem, Perdomo \cite{perdomo} achieved that assuming the hypersurfaces are preserved by the antipodal map.

The fourth and final problem is a beautiful conjecture of White \cite{white2} which says that any metric on a $3$-sphere has five distinct minimal embedded tori and he proved this for small perturbations  of the round metric \cite{white}. An easier conjecture would be to say that any metric on a $3$-sphere has nine distinct minimal surfaces of genus either zero or one.

\subsection*{Some new questions} For the purpose of applications in Geometry and Topology, it is important that to estimate the topology of the minimal hypersurface given by the Min-max Theorem.

 For ambient $3$-manifolds, due to the combined work of Simon--Smith \cite{smyth},  De Lellis--Pellandini \cite{Delellis-pellandini}, and Ketover \cite{ketover}, it is  now known that, roughly speaking, if the Min-max technique is applied to  continuous one-parameter family of embedded surfaces of genus $g$, then the min-max minimal surface has at most genus $g$.   

In higher dimensions it is not so clear how to control the topology of the min-max minimal hypersurface by the same methods.

An alternative approach would be to try to characterize the topology of the min-max minimal hypersurface via its index.  Note that if the minimal hypersurfaces are produced via Min-max methods then one should expect some control on the index.

 For ambient $3$-manifolds, Ejiri--Micallef \cite{ejiri-micallef} showed that the  index of a minimal orientable surface is bounded from above  by a multiple of area plus genus and if the metric has positive Ricci curvature then it is known \cite{yau2} that index one orientable minimal surfaces have genus $3$ at most (a conjecture that I heard from Rick Schoen states that the genus  should be two at most).

For higher dimensions, Rick Schoen conjectured that  index one embedded orientable  compact minimal hypersurfaces in ambient manifolds with positive Ricci curvature have bounded first Betti number.  We conjecture that if the ambient manifold has positive Ricci curvature then an index k embedded orientable compact minimal hypersurface has first Betti number bounded by fixed multiple of $k$.  Savo \cite{savo} showed this holds on round spheres of any dimension.

Another direction of research would be to understand the $p$-widths $\omega_p(M)$ of $(M^{n+1},g)$ for $n\geq 2$. The sequence $\{\omega_p(M)\}_{p\in\N}$ can be thought of as a nonlinear spectrum for the manifold and we would expect it to be asymptotically related with the spectrum of the Laplacian.  Taking this perspective, many interesting questions arise.

 For instance,  does the nonlinear spectrum satisfy a Weyl Law? More precisely, can we find a universal constant $a(n)$ so that
$$
\lim_{p\to\infty}\omega_p(M)p^{-\frac {1}{n+1}}=a(n)({\rm vol}(M,g))^{\frac{n}{n+1}}?
$$
This question has been suggested by Gromov in \cite[Section 8]{gromov} and in \cite[Section 5.2]{gromov2}. A classical result of Uhlenbeck \cite{uhlenbeck} states that  generic metrics have simple eigenvalues. Likewise, we expect that for   generic  metrics at least, $\omega_p(M)$ is achieved by a multiplicity one minimal hypersurface with index $p.$  Can we say anything about how they look like?  For instance, are they becoming equidistributed in space? The proof of Gromov-Guth's Theorem suggest that. Do they behave like nodal sets of eigenfunctions? Is their first betti number proportional to $p$?

Nodal sets of eigenfunctions provide a natural upper bound for $\omega_p(M)$. Making this more precise,  denote by $\phi_0, \dots, \phi_p$  the first $(p+1)$-eigenfunctions for the Laplace operator, where $\phi_0$. Consider
\begin{align*}
\Phi_p&:\RP^p\rightarrow\mathcal Z_n(M;\Z_2),\\
\Phi_p&([a_0,\ldots,a_p])=\partial\{x\in M:a_0\phi_0(x)+\ldots+a_p\phi_p(x)<0\}.
\end{align*}
and set
$$\bar \omega_p(M)=\sup_{\theta\in \RP^p}\vol(\Phi_p(\theta))\geq \omega_p(M).$$
It seems a challenging question to determine how  close  to one $\frac{\bar \omega_p(M)}{\omega_p(M)}$ is getting as $p$ tends to infinity. If the quotient is bounded, that would  imply a conjecture of Yau regarding the asymptotic   growth of the volume of nodal sets (which was proven in the analytic case  by Donnelly and Fefferman \cite{donnelly-fefferman} and for recent progress see  \cite{colding-minicozzi,sogge-zelditch}). Can we determine that quotient on an $n$-sphere?

\bibliographystyle{amsbook}

\end{document}